\documentclass[12pt]{article}%
\usepackage{amsmath}
\usepackage{amsfonts}
\usepackage{amssymb}
\usepackage{graphicx}%
\newtheorem{theorem}{Theorem}
\newtheorem{acknowledgement}[theorem]{Acknowledgement}
\newtheorem{algorithm}[theorem]{Algorithm}

\newtheorem{conclusion}[theorem]{Conclusion}

\newtheorem{corollary}[theorem]{Corollary}

\newtheorem{definition}[theorem]{Definition}
\newtheorem{example}[theorem]{Example}

\newtheorem{proposition}[theorem]{Proposition}
\newtheorem{remark}[theorem]{Remark}

\numberwithin{equation}{section}
\numberwithin{theorem}{section}
\newenvironment{proof}[1][Proof]{\noindent\textbf{#1.} }{\ \rule{0.5em}{0.5em}}

\def\lfhook#1{\setbox0=\hbox{#1}{\ooalign{\hidewidth
\lower1.5ex\hbox{'}\hidewidth\crcr\unhbox0}}}
\begin{document}

\title{Nested Derivatives: A simple method for computing series expansions of inverse functions.}
\author{Diego Dominici\thanks{dominicd@newpaltz.edu}\\Department of Mathematics\\State University of New York at New Paltz\\75 S. Manheim Blvd. Suite 9 \\New Paltz, NY 12561-2443, USA}
\maketitle

\begin{abstract}
We give an algorithm to compute the series expansion for the inverse of a
given function. The algorithm is extremely easy to implement and gives the
first $N$ terms of the series. We show several examples of its application in
calculating the inverses of some special functions.

\end{abstract}

MSC-class: 30B10 (Primary) 30D10, 33E20 (Secondary)

\section{Introduction}

\begin{quotation}
\textquotedblleft One must \textbf{always} invert.\textquotedblright

Carl G. J. Jacobi
\end{quotation}

\medskip

The existence of series expansions for inverses of analytic functions is a
well-known result of complex analysis \cite{MR31:2374}. \ The standard inverse
function theorem, a proof of which can be found, for example, in
\cite{MR86b:30001}, states that

\begin{theorem}
\label{theo 1}Let $h(x)$ be analytic for $\left\vert x-x_{0}\right\vert <R$
where $h^{\prime}(x_{0})\neq0.$ Then $z=h(x)$ has an analytic inverse $x=H(z)$
in some $\varepsilon$-neighborhood of $z_{0}=h(x_{0}).$
\end{theorem}

In the case when \ $x_{0}=z_{0}=0,$ $\left\vert h(x)\right\vert \leq M$ \ for
\ $\left\vert x\right\vert <R$, $\ $and $\ h^{\prime}(0)=a,$ \ R. M. Redheffer
\cite{redheffer} has shown that it is enough to take \ $\varepsilon=\frac
{1}{4}\frac{\left(  aR\right)  ^{2}}{M}.$ \ 

However, the \emph{procedure} to obtain the actual series is usually very
difficult to implement in practice. Under the conditions of Theorem
\ref{theo 1}, the two standard methods to compute the coefficients $b_{n}$ of
\[
h^{-1}(z)\equiv H(z)=\sum\limits_{n\geq0}b_{n}(z-z_{0})^{n}%
\]
are reversion of series \cite{MR23:B561}, \cite{reynolds}, \cite{orstrand},
and Lagrange's theorem. The first one requires to expand $h(x)$ around $x_{0}$%
\[
h(x)=\sum\limits_{n\geq0}a_{n}(x-x_{0})^{n}%
\]
and then solve for $b_{n}$ in the equation
\[
z=\sum\limits_{n\geq0}a_{n}\left[  \sum\limits_{n\geq0}b_{n}(z-z_{0}%
)^{n}-x_{0}\right]  ^{n}%
\]
by equating powers of $z$ and taking into account that \ $a_{0}=z_{0}$ \ and
\ $b_{0}=x_{0}$. This method is especially useful if all that is known about
$h(x)$ are the first few $a_{n}$. When \ $x_{0}=z_{0}=0$ and $a_{1}=a,$ it was
shown by E. T. Whittaker \cite{MR14:40c} that
\begin{align*}
b_{1} &  =\frac{1}{a},\quad b_{2}=-\frac{a_{2}}{a^{3}},\quad b_{3}=\frac
{1}{3!a^{5}}\left\vert
\begin{array}
[c]{cc}%
3a_{2} & a\\
6a_{3} & 4a_{2}%
\end{array}
\right\vert ,\quad\ldots\\
b_{n} &  =\frac{\left(  -1\right)  ^{n-1}}{n!a^{2n-1}}\left\vert
\begin{array}
[c]{ccccc}%
na_{2} & a & 0 & 0 & \cdots\\
2na_{3} & (n+1)a_{2} & 2a & 0 & \cdots\\
3na_{4} & (2n+1)a_{3} & (n+2)a_{2} & 3a & \cdots\\
4na_{5} & (3n+1)a_{4} & 2(n+1)a_{3} & (n+3)a_{2} & \cdots\\
\vdots & \vdots & \vdots & \vdots & \ddots
\end{array}
\right\vert
\end{align*}
where $\left\vert \cdot\right\vert \equiv\det\left(  \cdot\right)  .$ In
Example 10, we show how to get the $b_{n}$ in term of the $a_{n}$ using our method.

A computer system like Maple can reverse the power series of \ $h(x)$,
provided $h(x)$ is not too complicated, by using the command
\begin{align*}
&  >\text{Order}:=N+1;\\
&  >\text{solve ( series ( }h(x)\text{ , }x=x_{0}\text{ , }N+1\text{) =
}z\text{ , }x\text{);}%
\end{align*}
where $N$ is the number of terms wanted. Fast algorithms of order \ $(n\log
n)^{3/2}$ \ for reversion of series have been analyzed by Brent and Kung
\cite{MR58:25090}, \cite{MR55:1699}. The multivariate case has been studied by
several authors \cite{MR29:34}, \cite{MR97h:13018}, \cite{MR46:2824},
\cite{MR88h:05011} and Wright \cite{MR90d:13008} has studied the connection
between reversion of power series and \textquotedblleft rooted
trees\textquotedblright.

The second and more direct method is Lagrange's inversion formula
\cite{MR94b:00012},
\begin{equation}
b_{n}=\frac{1}{n!}\left.  \frac{d^{n-1}}{dx^{n-1}}\left\{  \left[
\frac{x-x_{0}}{h(x)-z_{0}}\right]  ^{n}\right\}  \right\vert _{x=x_{0}}.
\label{Lag}%
\end{equation}
Unfortunately, more direct doesn't necessarily mean easier and, except for
some simple cases, Lagrange's formula (\ref{Lag}) is extremely complicated for
practical applications. The $q$-analog (a mathematical expression parametrized
by $q$ which generalizes an expression and reduces to it in the limit
$q\rightarrow1^{+})$ of (\ref{Lag}) has been studied by various authors
\cite{MR52:10441}, \cite{MR82g:05006}, \cite{MR84f:33009}, \cite{MR87i:33007}
and a unified approach to both the regular and $q$-analog formulas have been
obtained by Krattenthaler \cite{MR89d:05017}. There has also been a great deal
of attention to the asymptotic expansion of inverses \cite{MR95e:11135},
\cite{MR2000h:41039}, \cite{MR93a:33003}, \cite{MR94e:33001}.

In this note, we present a simple, easy to implement method for computing the
series expansion for the inverse of any function satisfying the conditions of
Theorem \ref{theo 1}, although the method is especially powerful when $h(x)$
has the form
\[
h(x)=\int\limits_{a}^{x}g(x)dx
\]
and $g(x)$ is some function simpler than $h(x)$. Since this is the case for
many special functions, we will present several such examples. This note is
organized as follows:

In section 2 we define a sequence of functions $\mathfrak{D}^{n}[f]\,(x),$
obtained from a given one $f(x)$, that we call \textquotedblleft nested
derivatives\textquotedblright, for reasons which will be clear from the
definition. We give a computer code for generating the nested derivatives and
examples of how $\mathfrak{D}^{n}[f]\,(x)$ look for some elementary functions.
Section 3 shows how to compute the nested derivatives by using generating
functions. We present some examples and compare the results with those
obtained in Section 1.

Section 4 contains our main result of the use of nested derivatives to compute
power series of inverses. We test our result with some known results and we
apply the method for obtaining expansions for the inverse of the error
function, the incomplete Gamma function, the sine integral, and other special functions.

\section{Definitions}

\begin{definition}
We define $\mathfrak{D}^{n}[f]$\thinspace$(x),$ the n$^{th}$ nested derivative
of the function $f(x),$ by the following recursion:
\begin{align}
\mathfrak{D}^{0}[f]\,(x)  &  =1\nonumber\\
\mathfrak{D}^{n}[f]\,(x)  &  =\frac{d}{dx}\left[  f(x)\times\mathfrak{D}%
^{n-1}[f]\,(x)\right]  ,\quad n\geq1. \label{def}%
\end{align}

\end{definition}

\begin{proposition}
The nested derivative $\mathfrak{D}^{n}[f]$\thinspace$(x)$ satisfies the
following basic properties.

\begin{enumerate}
\item[(1)] For $n\geq1,\ \mathfrak{D}^{n}[\kappa]\equiv0,\quad\kappa$ constant.

\item[(2)] For $n\geq0,\ \mathfrak{D}^{n}[\kappa f]$\thinspace$(x)=\kappa
^{n}\mathfrak{D}^{n}[f]$\thinspace$(x),\quad\kappa$ constant.

\item[(3)] For $n\geq1,\ \mathfrak{D}^{n}[f]$\thinspace$(x)$ has the following
integral representation:
\[
\mathfrak{D}^{n}[f]\,(x)=\frac{1}{(2\pi i)^{n}}\oint\limits_{C_{1}}%
\oint\limits_{C_{2}}\cdots\oint\limits_{C_{n}}\frac{f(z_{n})}{(z_{n}-x)^{2}%
}\prod\limits_{k=1}^{n-1}\frac{f(z_{k})}{(z_{k}-z_{k+1})^{2}}dz_{n}\ldots
dz_{1},
\]
where $C_{k}$ is a small loop around $x$ in the complex plane.
\end{enumerate}
\end{proposition}

\begin{proof}
Properties (1) and (2) follow immediately from the definition of
$\mathfrak{D}^{n}[f]\,(x).$

To prove (3) we use induction on $n.$ For $n=1$ the result follows from
Cauchy's formula
\[
\mathfrak{D}^{1}[f]\,(x)\equiv\frac{df}{dx}=\frac{1}{2\pi i}\oint
\limits_{C_{1}}\frac{f(z_{1})}{(z_{1}-x)^{2}}dz_{1}.
\]
Assuming that the result is true for $n$ and using (\ref{def}), it follows that%

\begin{align*}
\mathfrak{D}^{n+1}[f]\,(x)  &  \equiv\frac{d}{dx}\left[  f(x)\times
\mathfrak{D}^{n}[f]\,(x)\right]  =\frac{1}{2\pi i}\oint\limits_{C_{n+1}}%
\frac{f(z_{n+1})\mathfrak{D}^{n}[f]\,(z_{n+1})}{(z_{n+1}-x)^{2}}dz_{n+1}\\
&  =\frac{1}{(2\pi i)^{n+1}}\oint\limits_{C_{1}}\cdots\oint\limits_{C_{n+1}%
}\frac{f(z_{n+1})}{(z_{n+1}-x)^{2}}\frac{f(z_{n})}{(z_{n}-z_{n+1})^{2}}\\
&  \times\prod\limits_{k=1}^{n-1}\frac{f(z_{k})}{(z_{k}-z_{k+1})^{2}}%
dz_{n+1}\ldots dz_{1}\\
&  =\frac{1}{(2\pi i)^{n+1}}\oint\limits_{C_{1}}\cdots\oint\limits_{C_{n+1}%
}\frac{f(z_{n+1})}{(z_{n+1}-x)^{2}}\prod\limits_{k=1}^{n}\frac{f(z_{k}%
)}{(z_{k}-z_{k+1})^{2}}dz_{n+1}\ldots dz_{1}.
\end{align*}

\end{proof}

\begin{algorithm}
The $\mathfrak{D}$ algorithm. The following Maple procedure implements the
recurrence relation (\ref{def}). We define $d(k)=\mathfrak{D}^{k}[f]\,(x),$
where $N$ is the number of terms desired.
\begin{align}
&  >d(0):=1\text{;}\nonumber\\
&  >\text{\texttt{for }}k\text{ \ \texttt{from} \ }0\text{ \texttt{to}
\ }N\text{ \texttt{do}:}\nonumber\\
&  >d(k+1)\text{: = \texttt{simplify} ( \texttt{diff} ( }f(x)\text{ }%
\ast\text{ }d(k)\text{, }x\text{ )):}\label{alg}\\
&  >\text{\texttt{print} ( }k+1\text{, }d(k+1)\text{ ):}\nonumber\\
&  >\text{\texttt{od}:}\nonumber
\end{align}

\end{algorithm}

\begin{example}
\label{ex f=x}The function $f(x)=x.$%

\begin{align*}
\mathfrak{D}^{1}[f]\,(x)  &  =1\\
\mathfrak{D}^{2}[f]\,(x)  &  =1\\
&  \vdots\\
\mathfrak{D}^{n}[f]\,(x)  &  =1.
\end{align*}

\end{example}

\begin{example}
\label{ex f=x^r}The power function $f(x)=x^{r},\quad r\neq1.$%
\begin{align*}
\mathfrak{D}^{1}[f]\,(x)  &  =rx^{r-1}\\
\mathfrak{D}^{2}[f]\,(x)  &  =r(2r-1)x^{2(r-1)}\\
\mathfrak{D}^{3}[f]\,(x)  &  =r(2r-1)(3r-2)x^{3(r-1)}\\
&  \vdots\\
\mathfrak{D}^{n}[f]\,(x)  &  =\prod\limits_{j=1}^{n}\left[  jr-(j-1)\right]
x^{n(r-1)}\\
&  =(r-1)^{n}\frac{\Gamma\left(  n+1+\frac{1}{r-1}\right)  }{\Gamma\left(
1+\frac{1}{r-1}\right)  }x^{n(r-1)}.
\end{align*}

Notice that when $r=\frac{k}{k+1},\quad k=1,2,\ldots,$ the sequence of nested
derivatives has only $k+1$ non-zero terms
\[
\mathfrak{D}^{n}[f]\,(x)=\left\{
\begin{array}
[c]{c}%
\frac{k!}{(k-n)!\ (k+1)^{n}}x^{-\frac{n}{k+1}},\quad0\leq n\leq k\\
0,\quad n\geq k+1
\end{array}
\right.  .
\]

\end{example}

\begin{example}
\label{ex f = exp}The exponenetial function $f(x)=e^{rx}.$%
\begin{align*}
\mathfrak{D}^{1}[f]\,(x)  &  =re^{rx}\\
\mathfrak{D}^{2}[f]\,(x)  &  =2r^{2}e^{2rx}\\
\mathfrak{D}^{3}[f]\,(x)  &  =6r^{3}e^{3rx}\\
&  \vdots\\
\mathfrak{D}^{n}[f]\,(x)  &  =n!r^{n}e^{nrx}.
\end{align*}

\end{example}

\section{Generating functions}

Generating functions provide a valuable method for computing sequences of
functions defined by an iterative process; we will use them to calculate
$\mathfrak{D}^{n}[f]\,(x).$ In the sequel, we shall implicitly assume that the
generating function series converges in some small disc around $z=0.$

\begin{theorem}
Given $h(x)=\int\frac{1}{f(x)}dx,$ its inverse$\ H(x)=h^{-1}(x)$\ and the
exponential generating function $G(x,z)=\sum\limits_{n\geq0}$ $\mathfrak{D}%
^{n}[f]\,(x)\frac{z^{n}}{n!}$, it follows that \
\begin{equation}
G(x,z)=\frac{1}{f(x)}(f\circ H)\left[  z+h(x)\right]  . \label{gener}%
\end{equation}

\end{theorem}

\begin{proof}
Taking (\ref{def}) into account gives%
\begin{align*}
\frac{\partial}{\partial x}\left[  f(x)\times G(x,z)\right]   &
=\sum\limits_{n\geq0}\frac{d}{dx}\left[  f(x)\times\mathfrak{D}^{n}%
[f]\,(x)\right]  \frac{z^{n}}{n!}\\
&  =\sum\limits_{n\geq0}\mathfrak{D}^{n+1}[f]\,(x)\frac{z^{n}}{n!}%
=\sum\limits_{n\geq1}\mathfrak{D}^{n}[f]\,(x)\frac{z^{n-1}}{(n-1)!}\\
&  =\frac{\partial}{\partial z}\sum\limits_{n\geq0}\mathfrak{D}^{n}%
[f]\,(x)\frac{z^{n}}{n!}=\frac{\partial}{\partial z}G(x,z).
\end{align*}

Hence, the generating function $G(x,z)$ satisfies the PDE
\[
\frac{\partial(f\times G)}{\partial x}=\frac{\partial G}{\partial z}%
\]
with general solution
\begin{equation}
G(x,z)=\frac{1}{f(x)}g\left[  z+h(x)\right]  \label{G1}%
\end{equation}
where $g(z)$ is an arbitrary analytic function. Invoking the boundary
condition $G(x,0)=\mathfrak{D}^{0}[f]\,(x)=1$, (\ref{G1}) gives
\[
1=\frac{1}{f(x)}g\left[  h(x)\right]
\]
and therefore
\[
f(x)=\left(  g\circ h\right)  (x).
\]
If $x=H(w),$ then
\[
\left(  f\circ H\right)  (w)=\left(  g\circ h\circ H\right)  (w)=g(w)
\]
and the theorem follows.
\end{proof}

\begin{example}
The function $f(x)=x.$

Here $h(x)=\int\frac{1}{x}dx=\ln(x),\quad H(x)=e^{x},$ and from (\ref{gener})
it follows that
\[
G(x,z)=\frac{1}{x}\exp\left[  z+\ln(x)\right]  =e^{z}.
\]

We could obtain the same result from Example \ref{ex f=x} by summing the
series
\[
G(x,z)=\sum\limits_{n\geq0}1\frac{z^{n}}{n!}=e^{z}.
\]

\end{example}

\begin{example}
The power function $f(x)=x^{r},\quad r\neq1.$ \ 

Now $h(x)=\int x^{-r}dx=\frac{x^{1-r}}{1-r},\quad H(x)=\left[  (1-r)x\right]
^{\frac{1}{1-r}},$ and we get
\begin{align*}
G(x,z)  &  =x^{-r}\left\{  \left[  (1-r)\left(  z+\frac{x^{1-r}}{1-r}\right)
\right]  ^{\frac{1}{1-r}}\right\}  ^{r}=\left[  \frac{(1-r)z+x^{1-r}}{x^{1-r}%
}\right]  ^{\frac{r}{1-r}}\\
&  =\left[  1+(1-r)x^{r-1}z\right]  ^{\frac{r}{1-r}}.
\end{align*}
Expanding in series around $z=0,$ we recover the result from Example
\ref{ex f=x^r}.

If $\frac{r}{1-r}=k,$ \ i.e. $r=\frac{k}{k+1},\quad k=0,1,\ldots$, \ then
$G(x,z)$ is a polynomial of degree $k$ in $z$ and hence
\[
\mathfrak{D}^{n}[f]\,(x)=0,\quad n\geq k+1
\]
as we have already observed in Example \ref{ex f=x^r}.
\end{example}

Given the particular form of the function $h(x)$ in Theorem 2, we can get
alternative expressions for (\ref{gener}) which sometimes are easier to employ.

\begin{corollary}
Let $h(x)=\int\frac{1}{f(x)}dx,$ its inverse $H(x)=h^{-1}(x)$ and the
exponential generating function $G(x,z)=\sum\limits_{n\geq0}$ $\mathfrak{D}%
^{n}[f]\,(x)\frac{z^{n}}{n!}.$ Then,

\begin{enumerate}
\item[(\textit{i})]
\begin{equation}
G(x,z)=\frac{1}{f(x)}H^{\prime}\left[  z+h(x)\right]  \label{corol1}%
\end{equation}
and

\item[(\textit{ii})]
\[
G(x,z)=\frac{d}{dx}H\left[  z+h(x)\right]  .
\]

\end{enumerate}
\end{corollary}

\begin{proof}

\begin{enumerate}
\item[(i)] By definition $\ (h\circ H)(x)=x$, so
\[
h^{\prime}\left[  H(x)\right]  H^{\prime}(x)=1.
\]
Since $h(x)=\int\frac{1}{f(t)}dt$,
\[
\frac{1}{f\left[  H(x)\right]  }H^{\prime}(x)=1
\]
or
\[
(f\circ H)(x)=H^{\prime}(x)
\]
and therefore
\begin{align*}
G(x,z)  &  =\frac{1}{f(x)}(f\circ H)\left[  z+h(x)\right] \\
&  =\frac{1}{f(x)}H^{\prime}\left[  z+h(x)\right]  .
\end{align*}

\item[(ii)] Using the chain rule
\begin{align*}
\frac{d}{dx}H\left[  z+h(x)\right]   &  =H^{\prime}\left[  z+h(x)\right]
h^{\prime}(x)\\
&  =H^{\prime}\left[  z+h(x)\right]  \frac{1}{f(x)}%
\end{align*}
and the conclusion follows from part (i).
\end{enumerate}
\end{proof}

\section{Applications}

We now state our main result.

\begin{theorem}
Let $h(x)=\int\limits_{a}^{x}\frac{1}{f(t)}dt,$ with $f(a)\neq0,\pm\infty
,\ $\ and its inverse $\ H(x)=h^{-1}(x).\ $Then,
\begin{equation}
H(z)=a+f(a)\sum\limits_{n\geq1}\mathfrak{D}^{n-1}[f]\,(a)\frac{z^{n}}{n!}
\label{Theo2}%
\end{equation}
where $\left\vert z\right\vert <\varepsilon,\ $for some $\varepsilon>0.$
\end{theorem}

\begin{proof}
Let's first observe that since $h(a)=0$, it follows that $H(0)=a$ and from
(\ref{corol1})
\[
G(a,z)=\frac{1}{f(a)}H^{\prime}\left[  z+h(a)\right]  =\frac{1}{f(a)}%
H^{\prime}(z)
\]
where
\[
G(a,z)=\sum\limits_{n\geq0}\mathfrak{D}^{n}[f]\,(a)\frac{z^{n}}{n!}.
\]
Hence,
\begin{align*}
H(z)  &  =H(0)+\int\limits_{0}^{z}f(a)\sum\limits_{n\geq0}\mathfrak{D}%
^{n}[f]\,(a)\frac{t^{n}}{n!}dt\\
&  =a+f(a)\sum\limits_{n\geq0}\mathfrak{D}^{n}[f]\,(a)\frac{z^{n+1}}{(n+1)!}\\
&  =a+f(a)\sum\limits_{n\geq1}\mathfrak{D}^{n-1}[f]\,(a)\frac{z^{n}}{n!}.
\end{align*}

\end{proof}

\begin{example}
The natural logarithm function. Let $f(x)=e^{-x},$ with $a=0.$ We have
\ $f(0)=1,$%
\[
h(x)=\int\limits_{0}^{x}e^{t}dt=e^{x}-1,\quad H(x)=\ln(x+1)
\]
and from Example \ref{ex f = exp}
\[
\mathfrak{D}^{n}[f]\,(0)=(-1)^{n}n!.
\]

Hence, from (\ref{Theo2}) we get the familiar formula
\[
\ln(z+1)=\sum\limits_{n\geq1}(-1)^{n-1}\frac{z^{n}}{n}.
\]

\end{example}

\begin{example}
The tangent function. Let $f(x)=x^{2}+1,$ with $a=0.$Now $f(0)=1,$%
\[
h(x)=\int\limits_{0}^{x}\frac{1}{t^{2}+1}dt=\arctan(x),\quad H(x)=\tan(x)
\]
and (\ref{Theo2}) implies that
\[
\tan(z)=\sum\limits_{n\geq1}\mathfrak{D}^{n-1}[x^{2}+1]\,(0)\frac{z^{n}}{n!}.
\]

Therefore,
\begin{align}
\mathfrak{D}^{2k+1}[x^{2}+1]\,(0)  &  =0,\quad k\geq0\nonumber\\
\mathfrak{D}^{2k}[x^{2}+1]\,(0)  &  =\frac{2}{k+1}4^{k}\left(  4^{k+1}%
-1\right)  \left\vert B_{2(k+1)}\right\vert ,\quad k\geq1 \label{x^2+1}%
\end{align}
where $B_{k}$ are the Bernoulli numbers \cite{MR94b:00012}.
\end{example}

\begin{remark}
From Example \ref{ex f=x^r}, we recall that
\[
\mathfrak{D}^{n}[x^{2}]\,(x)=(n+1)!x^{n}%
\]
and consequently
\begin{equation}
\mathfrak{D}^{n}[x^{2}]\,(0)=0,\quad n\geq1. \label{x^2}%
\end{equation}

Comparing (\ref{x^2+1}) and (\ref{x^2}) we can see the highly nonlinear
behavior of \ the nested derivatives, since even the addition of $\ 1$ to
$f(x)$ creates a completely different sequence of values, far more complex
than the original.
\end{remark}

We now start testing our result on some classical functions.

\begin{example}
Elliptic functions. Let $f(x)=\sqrt{1-p^{2}\sin^{2}(x)},\quad0\leq
p\leq1,\quad a=0.$ We have, \ $f(0)=1$ \ and
\[
h(\phi)=\int\limits_{0}^{\phi}\frac{d\theta}{\sqrt{1-p^{2}\sin^{2}(\theta)}%
}=F(p;\phi),\quad H(p;x)=\operatorname{am}(p;x)
\]
where $F(p;\phi)$ is the incomplete elliptic integral of the first kind, and
$\operatorname{am}(p;x)$ is the elliptic amplitude \cite{Atlas}
\[
\operatorname{am}(p;x)=\arcsin\left[  \operatorname{sn}(p;x)\right]
=\arccos\left[  cn(p;x)\right]  =\arcsin\left[  \frac{\sqrt{1-dn^{2}(p;x)}}%
{p}\right]
\]
with $\operatorname{sn}(p;x),\ cn(p;x),\ $and $dn(p;x)$ denoting the Jacobian
elliptic functions.

Computing $\mathfrak{D}^{n}[f]\,(0)$ with (\ref{alg}) gives
\begin{align*}
\mathfrak{D}^{2k+1}[f]\,(0)  &  =0,\quad k\geq0\\
\mathfrak{D}^{2k}[f]\,(0)  &  =(-1)^{k}p^{2}Q_{k}(p),\quad k\geq1
\end{align*}
where $Q_{k}(p)$ is a polynomial of degree $2(k-1)$ of the form
\[
Q_{k}(p)=p^{2(k-1)}+\cdots+2^{2(k-1)}.
\]
The first few $Q_{k}(p)$ \ are
\begin{align*}
Q_{1}(p)  &  =1\\
Q_{2}(p)  &  =p^{2}+4\\
Q_{3}(p)  &  =p^{4}+44p^{2}+16\\
Q_{4}(p)  &  =p^{6}+408p^{4}+912p^{2}+64\\
Q_{5}(p)  &  =p^{8}+3688p^{6}+307682p^{4}+15808p^{2}+256
\end{align*}
and (\ref{Theo2}) implies that
\begin{equation}
am(p;x)=z-p^{2}\frac{z^{3}}{3!}+p^{2}(p^{2}+4)\frac{z^{5}}{5!}-p^{2}%
(p^{4}+44p^{2}+16)\frac{z^{7}}{7!}+\cdots\label{am}%
\end{equation}
in agreement with the known expansions for $am(p;x)$ \cite{MR31:6000}$.$
\end{example}

\begin{example}
The Lambert-W function. Let $f(x)=e^{-x}(x+1)^{-1},\quad a=0,\quad f(a)=1.$
Here
\[
h(x)=xe^{x},\quad H(x)=LW(x)
\]
where by $LW(x)$ we denote the Lambert-W function \cite{MR98j:33015},
\cite{MR1809988}, \cite{MR97e:33003}. In this case, (\ref{alg}) gives
\begin{align*}
\mathfrak{D}^{1}[f]\,(0)  &  =-2\\
\mathfrak{D}^{2}[f]\,(0)  &  =9\\
\mathfrak{D}^{3}[f]\,(0)  &  =-64\\
&  \vdots\\
\mathfrak{D}^{n}[f]\,(0)  &  =\left[  -(n+1)\right]  ^{n}.
\end{align*}

From (\ref{Theo2}) we conclude that
\[
LW(z)=\sum\limits_{n\geq1}(-1)^{n-1}n^{n-1}\frac{z^{n}}{n!}.
\]

\end{example}

\begin{example}
We now derive a well known result \cite{MR94b:00012} about reversion of
series. If we take
\[
h(x)=a_{1}x+a_{2}x^{2}+a_{3}x^{3}+a_{4}x^{4}+a_{5}x^{5}+a_{6}x^{6}+a_{7}%
x^{7}+\cdots
\]
where $a_{1}\neq0,\ $then
\[
f(x)=\frac{1}{h^{\prime}(x)}=\frac{1}{a_{1}+2a_{2}x+3a_{3}x^{2}+4a_{4}%
x^{3}+5a_{5}x^{4}+6a_{6}x^{5}+\allowbreak7a_{7}x^{6}+\cdots}%
\]
$a=0,\ f(0)=\frac{1}{a_{1}}\ $and from (\ref{alg}) we get
\begin{align*}
\mathfrak{D}^{1}[f]\,(0)  &  =-2\frac{a_{2}}{\left(  a_{1}\right)  ^{2}}\\
\mathfrak{D}^{2}[f]\,(0)  &  =6\frac{2\left(  a_{2}\right)  ^{2}-a_{1}a_{3}%
}{\left(  a_{1}\right)  ^{4}}\\
\mathfrak{D}^{3}[f]\,(0)  &  =24\frac{5a_{1}a_{2}a_{3}-\left(  a_{1}\right)
^{2}a_{4}-5(a_{2})^{3}}{\left(  a_{1}\right)  ^{6}}\\
\mathfrak{D}^{4}[f]\,(0)  &  =120\frac{6(a_{1})^{2}a_{2}a_{3}+3\left(
a_{1}a_{3}\right)  ^{2}+14(a_{2})^{4}-(a_{1})^{3}a_{5}-21a_{1}(a_{2})^{2}%
a_{3}}{\left(  a_{1}\right)  ^{8}}.
\end{align*}

Hence,
\begin{align*}
H(z)  &  =\frac{1}{a_{1}}z-\frac{a_{2}}{\left(  a_{1}\right)  ^{3}}z^{2}%
+\frac{2\left(  a_{2}\right)  ^{2}-a_{1}a_{3}}{\left(  a_{1}\right)  ^{5}%
}z^{3}+\frac{5a_{1}a_{2}a_{3}-\left(  a_{1}\right)  ^{2}a_{4}-5(a_{2})^{3}%
}{\left(  a_{1}\right)  ^{7}}z^{4}\\
&  +\frac{6(a_{1})^{2}a_{2}a_{3}+3\left(  a_{1}a_{3}\right)  ^{2}%
+14(a_{2})^{4}-(a_{1})^{3}a_{5}-21a_{1}(a_{2})^{2}a_{3}}{\left(  a_{1}\right)
^{9}}z^{5}+\cdots.
\end{align*}

\end{example}

\begin{remark}
An explicit formula for the $n^{\text{th}}$ term is given in Morse and
Feshbach \cite[Part 1 pp. 411--413]{MR15:583h},%
\begin{gather*}
b_{n}=\frac{1}{n\left(  a_{1}\right)  ^{n}}%
{\displaystyle\sum\limits_{s,t,u,\ldots}}
\left(  -1\right)  ^{s+t+u+\cdots}\frac{n(n+1)\cdots(n-1+s+t+u+\cdots
)}{s!t!u!\cdots}\left(  \frac{a_{2}}{a_{1}}\right)  ^{s}\left(  \frac{a_{3}%
}{a_{1}}\right)  ^{t}\cdots\\
s+2t+3u+\cdots=n-1.
\end{gather*}

\end{remark}

\begin{example}
The Error Function, $\operatorname{erf}(x).$ We now have
\begin{align*}
h(x)  &  =\operatorname{erf}(x)=\frac{2}{\sqrt{\pi}}\int\limits_{0}%
^{x}e^{-t^{2}}dt\\
f(x)  &  =\frac{\sqrt{\pi}}{2}e^{x^{2}},\quad a=0,\quad f(a)=\frac{\sqrt{\pi}%
}{2}%
\end{align*}
and (\ref{alg}) gives
\[
\mathfrak{D}^{n}[f]\,(0)=\left\{
\begin{array}
[c]{c}%
0,\quad n=2k+1,\quad k\geq0\\
\left(  \frac{\sqrt{\pi}}{2}\right)  ^{n}A_{k},\quad n=2k,\quad k\geq0
\end{array}
\right.
\]
where
\begin{align*}
A_{0}  &  =1,\quad A_{1}=2,\quad A_{2}=28,\quad A_{3}=1016,\quad A_{4}=69904\\
A_{5}  &  =7796768,\quad A_{6}=1282366912,\quad A_{7}=291885678464,\ldots.
\end{align*}

From (\ref{Theo2}) we get
\[
H(z)=\sum\limits_{n\geq0}A_{n}\left(  \frac{\sqrt{\pi}}{2}\right)
^{2n+1}\frac{z^{2n+1}}{(2n+1)!}.
\]

which agrees with other authors calculations previously published
\cite{MR54:9047}, \cite{MR27:3839}, \cite{Inverse}, \cite{MR49:6558},
\cite{MR36:6119}.
\end{example}

We will now extend (\ref{Theo2}) to a more general result.

\begin{corollary}
Let $h(x)=\int\limits_{a}^{x}\frac{1}{f(t)}dt,$ $\ z_{0}=h(b),\ $with
$f(b)\neq0,\pm\infty$\ and its inverse $\ H(x)=h^{-1}(x).$ Then,%
\begin{equation}
H(z)=b+f(b)\sum\limits_{n\geq1}\mathfrak{D}^{n-1}[f]\,(b)\frac{(z-z_{0})^{n}%
}{n!} \label{corol2}%
\end{equation}
where $\left\vert z-z_{0}\right\vert <\varepsilon,\ $for some $\varepsilon>0.$
\end{corollary}

\begin{proof}
We consider the function
\[
u(x)=h(x)-z_{0}%
\]
which satisfies $\ u(b)=0,\ $and its inverse $U(x)=u^{-1}(x).$ \ Since
$f(b)\neq0,\pm\infty,\ $we can apply (\ref{Theo2})$\ $to $u(x)$ and conclude
that
\[
U(z)=b+f(b)\sum\limits_{n\geq1}\mathfrak{D}^{n-1}[f]\,(b)\frac{z^{n}}{n!}.
\]

All that is left is to see the relation between $U(z)$ and $H(z).$ \ 

Suppose that $u(x)=y$. Then
\begin{align*}
y  &  =u(x)=h(x)-z_{0}\\
h(x)  &  =y+z_{0}\\
x  &  =H(y+z_{0})
\end{align*}
and therefore
\[
U(y)=H(y+z_{0})
\]
or
\[
H(z)=U(z-z_{0})
\]
and (\ref{corol2}) follows.
\end{proof}

We will now use our results to get some power series expansions that have not
been studied before.

\begin{example}
The incomplete Gamma function, $\gamma(\nu;x).$ We have
\begin{align*}
h(\nu;x)  &  =\gamma(\nu;x)\equiv\int\limits_{0}^{x}e^{-t}t^{\nu-1}dt,\quad
\nu>0,\quad x\geq0\\
f(\nu;x)  &  =e^{x}x^{1-\nu},\quad a=0.
\end{align*}
Since
\[
f(\nu;0)=\left\{
\begin{array}
[c]{c}%
0,\quad0<\nu<1\\
\infty,\quad\nu>1
\end{array}
\right.
\]
we can't apply (\ref{Theo2}). Choosing $b=1,\quad z_{0}(\nu)=\gamma
(\nu;1),\quad f(\nu;b)=e,$ we conclude from (\ref{corol2}) that%
\[
H(\nu;z)=1+e\sum\limits_{n\geq1}\mathfrak{D}^{n-1}[f]\,(1)\frac{\left[
z-z_{0}(\nu)\right]  ^{n}}{n!}.
\]

We use (\ref{alg}) to compute the first few $\mathfrak{D}^{n}[f]\,(1)$ and
obtain
\[
\mathfrak{D}^{n}[f]\,(1)=e^{n}Q_{n}(\nu)
\]
where $Q_{n}(\nu)$ is a polynomial of degree $n$%
\begin{align*}
Q_{1}(\nu)  &  =2-\nu\\
Q_{2}(\nu)  &  =7-7\nu+2\nu^{2}\\
Q_{3}(\nu)  &  =36-53\nu+29\nu^{2}-6\nu^{3}\\
Q_{4}(\nu)  &  =\allowbreak245-474\nu+375\nu^{2}-146\nu^{3}+24\nu
^{4}\allowbreak\\
Q_{5}(\nu)  &  =2076-4967\nu+5104\nu^{2}-2847\nu^{3}+874\nu^{4}-\allowbreak
120\nu^{5}\allowbreak
\end{align*}

and we can write
\[
H(\nu;z)=1+\sum\limits_{n\geq1}e^{n}Q_{n-1}(\nu)\frac{\left[  z-z_{0}%
(\nu)\right]  ^{n}}{n!}.
\]

\end{example}

\begin{example}
The sine integral function, $\operatorname{Si}(x).$ In this case
\[
h(x)=\operatorname{Si}(x)\equiv\int\limits_{0}^{x}\frac{\sin(t)}{t}dt,\quad
f(x)=\frac{x}{\sin(x)},\quad a=0.
\]

For this example $f(a)$ is well defined, but to simplify the calculations we
choose $b=\frac{\pi}{2},\quad z_{0}=\operatorname{Si}(\frac{\pi}{2})\simeq$
$1.370762.$ Then,
\[
f(b)=\frac{\pi}{2},\quad\mathfrak{D}^{n}[f]\,\left(  \frac{\pi}{2}\right)
=Q_{n}(\pi)
\]
where $Q_{n}(x)$ is once again a polynomial
\begin{align*}
Q_{1}(x)  &  =1\\
Q_{2}(x)  &  =1+\frac{1}{4}x^{2}\\
Q_{3}(x)  &  =1+\frac{7}{4}x^{2}\\
Q_{4}(x)  &  =1+8x^{2}+\frac{9}{16}x^{4}\\
Q_{5}(x)  &  =1+\frac{61}{2}x^{2}+\frac{159}{16}x^{4}\\
Q_{6}(x)  &  =1+\frac{423}{4}x^{2}+\frac{1671}{16}x^{4}+\frac{225}{64}x^{6}.
\end{align*}
and from (\ref{corol2}) we obtain
\[
H(z)=\frac{\pi}{2}+\frac{\pi}{2}\sum\limits_{n\geq1}Q_{n}(\pi)\frac
{(z-z_{0})^{n}}{n!}.
\]

\end{example}

\begin{example}
The logarithm integral function, $\operatorname{li}(x).$ From the definition
\[
h(x)=\mathrm{li}(x)\equiv\int\limits_{0}^{x}\frac{1}{\ln(t)}dt,\quad
f(x)=\ln(x),\quad a=0.
\]

In this case $f(a)=-\infty,$ so we must choose $b$. A natural candidate is
$b=e,$ which gives
\begin{align*}
f(b)  &  =1,\quad z_{0}=\mathrm{li}(e)\simeq1.895117816\\
\mathfrak{D}^{n}[f]\,(e)  &  =e^{-n}A_{n}%
\end{align*}
with
\begin{align*}
A_{1}  &  =1,\quad A_{2}=0,\quad A_{3}=-1,\quad A_{4}=2,\quad A_{5}=1\\
A_{6}  &  =-26,\quad A_{7}=99,\quad A_{8}=90,\quad A_{9}=-3627
\end{align*}
and we have
\[
H(z)=e+\sum\limits_{n\geq1}A_{n}\frac{(z-z_{0})^{n}}{n!}.
\]

\end{example}

\begin{example}
The incomplete Beta function, $B(\nu,\mu;x).$ By definition
\[
h(\nu,\mu;x)=B(\nu,\mu;x)\equiv\int\limits_{0}^{x}t^{\nu-1}(1-t)^{\mu
-1}dt,\quad0\leq x<1
\]
and hence
\[
f(\nu,\mu;x)=x^{1-\nu}(1-x)^{1-\mu},\quad a=0.
\]

To avoid the possible singularities at $x=0$ and $x=1$ we consider $b=\frac
{1}{2},$\ and therefore
\[
f(\nu,\mu;b)=\frac{1}{4}2^{\nu+\mu},\quad z_{0}(\nu,\mu)=B(\nu,\mu;\frac{1}%
{2}).
\]

The $\mathfrak{D}$ algorithm now gives
\[
\mathfrak{D}^{n}[f]\,\left(  \frac{1}{2}\right)  =2^{n(\nu+\mu-1)}Q_{n}%
(\nu,\mu)
\]
with $Q_{n}(\nu,\mu)$ a multivariate polynomial of degree $n$
\begin{align*}
Q_{1}(\nu,\mu)  &  =\mu-\nu\\
Q_{2}(\nu,\mu)  &  =-2+\nu-4\nu\mu+\mu+2\mu^{2}+2\nu^{2}\\
Q_{3}(\nu,\mu)  &  =(\mu-\nu)(6\mu^{2}-12\nu\mu+7\mu-12+6\nu^{2}+7\nu)\\
Q_{4}(\nu,\mu)  &  =16-46\nu\mu^{2}-46\nu^{2}\mu-63\mu^{2}-22\mu+154\nu
\mu-96\nu\mu^{3}-96\nu^{3}\mu+\\
&  144\nu^{2}\mu^{2}-22\nu-63\nu^{2}+24\nu^{4}+46\nu^{3}+24\mu^{4}+46\mu^{3}\\
Q_{5}(\nu,\mu)  &  =(\mu-\nu)(120\mu^{4}+326\mu^{3}-480\nu\mu^{3}+720\nu
^{2}\mu^{2}-323\mu^{2}-326\nu\mu^{2}\\
&  -362\mu-480\nu^{3}\mu+1154\nu\mu-326\nu^{2}\mu-323\nu^{2}+240+120\nu^{4}\\
&  +326\nu^{3}-362\nu).
\end{align*}
and we have
\[
H(z)=\frac{1}{2}+\frac{1}{4}2^{\nu+\mu}\sum\limits_{n\geq1}2^{n(\nu+\mu
-1)}Q_{n}(\nu,\mu)\frac{(z-z_{0})^{n}}{n!}.
\]

\end{example}

\begin{conclusion}
We have presented a simple method for computing the series expansion for the
inverses of functions and given a Maple procedure to generate the coefficients
in these expansions. We showed several examples of the method applied to
elementary and special functions, and stated the first few terms of the series
in each case.
\end{conclusion}

\begin{acknowledgement}
I wish to express my gratitude to the referees for their extremely valuable
suggestions and comments on previous versions of this work.
\end{acknowledgement}


\begin{thebibliography}{10}

\bibitem{MR94b:00012}
M.~Abramowitz and I.~A. Stegun, editors.
\newblock {\em Handbook of mathematical functions with formulas, graphs, and
  mathematical tables}.
\newblock Dover Publications Inc., New York, 1992.
\newblock Reprint of the 1972 edition.

\bibitem{MR52:10441}
G.~E. Andrews.
\newblock Identities in combinatorics. {II}. {A} {$q$}-analog of the {L}agrange
  inversion theorem.
\newblock {\em Proc. Amer. Math. Soc.}, 53(1):240--245, 1975.

\bibitem{MR54:9047}
J.~M. Blair, C.~A. Edwards, and J.~H. Johnson.
\newblock Rational {C}hebyshev approximations for the inverse of the error
  function.
\newblock {\em Math. Comp.}, 30(136):827--830, 1976.

\bibitem{MR29:34}
G.~R. Blakley.
\newblock Formal solution of nonlinear simultaneous equations: {R}eversion of
  series in several variables.
\newblock {\em Duke Math. J.}, 31:347--357, 1964.

\bibitem{MR55:1699}
R.~P. Brent and H.~T. Kung.
\newblock {$O((n$} {${\rm log} n)\sp{3/2})$} algorithms for composition and
  reversion of power series.
\newblock In {\em Analytic computational complexity (Proc. Sympos.,
  Carnegie-Mellon Univ., Pittsburgh, Pa., 1975)}, pages 217--225. Academic
  Press, New York, 1976.

\bibitem{MR58:25090}
R.~P. Brent and H.~T. Kung.
\newblock Fast algorithms for manipulating formal power series.
\newblock {\em J. Assoc. Comput. Mach.}, 25(4):581--595, 1978.

\bibitem{MR27:3839}
L.~Carlitz.
\newblock The inverse of the error function.
\newblock {\em Pacific J. Math.}, 13:459--470, 1963.

\bibitem{MR97h:13018}
C.~C.-A. Cheng, J.~H. McKay, J.~Towber, S.~S.-S. Wang, and D.~L. Wright.
\newblock Reversion of power series and the extended {R}aney coefficients.
\newblock {\em Trans. Amer. Math. Soc.}, 349(5):1769--1782, 1997.

\bibitem{MR98j:33015}
R.~M. Corless, G.~H. Gonnet, D.~E.~G. Hare, D.~J. Jeffrey, and D.~E. Knuth.
\newblock On the {L}ambert {$W$} function.
\newblock {\em Adv. Comput. Math.}, 5(4):329--359, 1996.

\bibitem{MR1809988}
R.~M. Corless, D.~J. Jeffrey, and D.~E. Knuth.
\newblock A sequence of series for the {L}ambert {$W$} function.
\newblock In {\em Proceedings of the 1997 International Symposium on Symbolic
  and Algebraic Computation (Kihei, HI)}, pages 197--204 (electronic), New
  York, 1997. ACM.

\bibitem{MR31:6000}
H.~T. Davis.
\newblock {\em Introduction to nonlinear differential and integral equations}.
\newblock Dover Publications Inc., New York, 1962.

\bibitem{MR86b:30001}
J.~W. Dettman.
\newblock {\em Applied complex variables}.
\newblock Dover Publications Inc., New York, 1984.
\newblock Reprint of the 1965 original.

\bibitem{Inverse}
D.~Dominici.
\newblock The inverse of the cumulative standard normal probability function.
\newblock {\em Integral Transforms Spec. Funct.}, 14(4):281--292, 2003.

\bibitem{MR46:2824}
R.~H. Estes and E.~R. Lancaster.
\newblock Some generalized power series inversions.
\newblock {\em SIAM J. Numer. Anal.}, 9:241--247, 1972.

\bibitem{MR49:6558}
H.~E. Fettis.
\newblock A stable algorithm for computing the inverse error function in the
  ``tail-end'' region.
\newblock {\em Math. Comp.}, 28:585--587, 1974.

\bibitem{MR23:B561}
D.~C. Fielder.
\newblock Tabulation of coefficients for operations on {T}aylor series.
\newblock {\em Math. Comp.}, 14:339--345, 1960.

\bibitem{MR31:2374}
A.~R. Forsyth.
\newblock {\em Theory of functions of a complex variable. {V}ols. 1 and 2}.
\newblock Third edition. Dover Publications Inc., New York, 1965.

\bibitem{MR82g:05006}
I.~Gessel.
\newblock A noncommutative generalization and {$q$}-analog of the {L}agrange
  inversion formula.
\newblock {\em Trans. Amer. Math. Soc.}, 257(2):455--482, 1980.

\bibitem{MR84f:33009}
I.~Gessel and D.~Stanton.
\newblock Applications of {$q$}-{L}agrange inversion to basic hypergeometric
  series.
\newblock {\em Trans. Amer. Math. Soc.}, 277(1):173--201, 1983.

\bibitem{MR87i:33007}
I.~Gessel and D.~Stanton.
\newblock Another family of {$q$}-{L}agrange inversion formulas.
\newblock {\em Rocky Mountain J. Math.}, 16(2):373--384, 1986.

\bibitem{MR88h:05011}
I.~M. Gessel.
\newblock A combinatorial proof of the multivariable {L}agrange inversion
  formula.
\newblock {\em J. Combin. Theory Ser. A}, 45(2):178--195, 1987.

\bibitem{MR97e:33003}
D.~J. Jeffrey, D.~E.~G. Hare, and R.~M. Corless.
\newblock Unwinding the branches of the {L}ambert {$W$} function.
\newblock {\em Math. Sci.}, 21(1):1--7, 1996.

\bibitem{MR89d:05017}
C.~Krattenthaler.
\newblock Operator methods and {L}agrange inversion: a unified approach to
  {L}agrange formulas.
\newblock {\em Trans. Amer. Math. Soc.}, 305(2):431--465, 1988.

\bibitem{MR15:583h}
P.~M. Morse and H.~Feshbach.
\newblock {\em Methods of theoretical physics. 2 volumes}.
\newblock McGraw-Hill Book Co., Inc., New York, 1953.

\bibitem{redheffer}
R.~M. Redheffer.
\newblock Reversion of power series.
\newblock {\em Amer. Math. Monthly}, 69(5):423--425, 1962.

\bibitem{reynolds}
J.~B. Reynolds.
\newblock Reversion of series with applications.
\newblock {\em Amer. Math. Monthly}, 51(10):578--580, 1944.

\bibitem{MR95e:11135}
B.~Salvy.
\newblock Fast computation of some asymptotic functional inverses.
\newblock {\em J. Symbolic Comput.}, 17(3):227--236, 1994.

\bibitem{MR2000h:41039}
B.~Salvy and J.~Shackell.
\newblock Symbolic asymptotics: multiseries of inverse functions.
\newblock {\em J. Symbolic Comput.}, 27(6):543--563, 1999.

\bibitem{Atlas}
J.~Spanier and K.~B. Oldham.
\newblock {\em An Atlas of Functions}.
\newblock Hemisphere Publishing Corp., New York, 1987.

\bibitem{MR36:6119}
A.~Strecok.
\newblock On the calculation of the inverse of the error function.
\newblock {\em Math. Comp.}, 22:144--158, 1968.

\bibitem{MR93a:33003}
N.~M. Temme.
\newblock Asymptotic inversion of incomplete gamma functions.
\newblock {\em Math. Comp.}, 58(198):755--764, 1992.

\bibitem{MR94e:33001}
N.~M. Temme.
\newblock Asymptotic inversion of the incomplete beta function.
\newblock {\em J. Comput. Appl. Math.}, 41(1-2):145--157, 1992.
\newblock Asymptotic methods in analysis and combinatorics.

\bibitem{orstrand}
C.~E. Van~Orstrand.
\newblock Reversion of power series.
\newblock 19(109):366--376, 1910.

\bibitem{MR14:40c}
E.~T. Whittaker.
\newblock On the reversion of series.
\newblock {\em Gaz. Mat., Lisboa}, 12(50):1, 1951.

\bibitem{MR90d:13008}
D.~Wright.
\newblock The tree formulas for reversion of power series.
\newblock {\em J. Pure Appl. Algebra}, 57(2):191--211, 1989.

\end{thebibliography}
\end{document}